\documentclass[12pt]{amsart}
\usepackage{amscd,amssymb,graphics}
\usepackage[matrix,arrow,curve]{xy}
\newcommand\mathscr{\mathcal}

\newtheorem{theorem}{Theorem}[section]
\newtheorem{corollary}[theorem]{Corollary} 

\newtheorem{proposition}[theorem]{Proposition}
\theoremstyle{definition}

\theoremstyle{remark}

\theoremstyle{fact}
\newtheorem{fact}[theorem]{Fact}
\newtheorem{conjecture}[theorem]{Conjecture}
\newtheorem{question}[theorem]{Question}

\renewcommand{\Bbb}{\mathbb}

\def\C {{\Bbb C}}

\def\QED{\nobreak\quad\ifmmode\roman{Q.E.D.}\else{\rm Q.E.D.}\fi}

\def\R{{\mathbf R}}

\def\o{\omega}

\def\sT{{\mathcal T}}

\def\sbs{\subset}

\def\ti{\times}
\def\obr{^{-1}}
\def\stm{\setminus}
\def\Card{\mathrm{Card}\,}
\def\Iso{\mathrm{Iso\,}}

\begin{document}

\title[Unitary representability]
{Unitary representability of free abelian topological groups}

\author[V.V. Uspenskij]{Vladimir V. Uspenskij}

\address{Department of Mathematics, 321 Morton Hall, Ohio
University, Athens, Ohio 45701, USA}

\email{uspensk@math.ohiou.edu}


\date{11 April 2006}


\begin{abstract} 
For every Tikhonov space $X$ the free abelian topological group
$A(X)$ and the free locally convex vector space $L(X)$ admit
a topologically faithful unitary representation. For 
compact spaces $X$ this is due to Jorge Galindo.
\end{abstract}

\maketitle

\section{Introduction} \label{s:intro}
With every Tikhonov space $X$ one associates the free topological group
$F(X)$, the free abelian topological group $A(X)$, and the free locally
convex vector space $L(X)$. They are characterized by respective universal
properties. For example, $L(X)$ is defined by the following: $X$ is an
(algebraic) basis of $L(X)$, and for every continuous mapping 
$f:X\to E$, where $E$ is a Hausdorff locally convex space, 
the  linear extension $\bar f:L(X)\to E$ of $f$ is continuous. There are two versions
of $L(X)$, real and complex. There also are versions of all these free
objects for spaces with a distinguished point. We consider non-pointed
spaces.

The {\em unitary group} $U(H)$ of a Hilbert space $H$ will be equipped
with the strong operator topology, which is the topology inherited
from the Tikhonov product $H^H$, or else the topology of pointwise
convergence. We use the notation $U_s(H)$ to indicate this topology. 
A {\em unitary representation} of a topological group $G$ is 
a continuous homomorphism $f:G\to U_s(H)$. Such a representation is
{\em faithful} if $f$ is injective, and {\em topologically faithful}
if $f$ is a homeomorphic embedding. A topological group is {\em unitarily
representable} if it is isomorphic to a topological subgroup of $U_s(H)$
for a Hilbert space $H$ (which may be non-separable), or, equivalently, 
if it admits a topologically faithful unitary representation.

All locally compact groups are unitarily representable.
For groups beyond the class of locally compact groups
this may no longer be true: there exist abelian topological groups
(even monothetic groups, that is, topologically generated by one element)
for which every unitary representation is trivial (sends the whole
group to the identity), see 
\cite[Theorem 5.1 and Remark 5.2]{Ban83}. 
Thus one may wonder what happens in the case
of free topological groups: are they unitarily representable?

In the non-abelian case, this question is open even if $X$
is compact metric, see \cite[Questions 35, 36]{P-40Q}. 
The aim of the present note is to answer the question
in the positive for $L(X)$ and $A(X)$.

\begin{theorem}
\label{th1}
For every Tikhonov space $X$ the free locally convex space $L(X)$ and
the free abelian topological group $A(X)$ admit a topologically faithful
unitary representation.
\end{theorem}

It suffices to consider the case of $L(X)$, since $A(X)$ is isomorphic 
to the subgroup of $L(X)$ generated by $X$ \cite{Tka}, \cite{Usp90}, see also
\cite{PGao02}
\footnote{The result was stated in \cite{Tka}, but the proof contained a gap.
I filled it in \cite{Usp90} --- not knowing that I was rediscovering the 
Monge -- Kantorovich metric (which is also known under many other names:
Prokhorov, Wasserstein, transportation, Earth Mover's) and its basic properties, 
such as the Integer Value
Property. The proof is reproduced in \cite{PGao02}, where a historic account
is given.}. For compact $X$ (or, more generally, for $k_\o$-spaces $X$)
Theorem~\ref{th1} is due to Jorge Galindo \cite{Gal05}. However, it was claimed
in \cite{Gal05} 
that there exists a metrizable space $X$ for which the group $A(X)$
is not isomorphic to a subgroup of a unitary group. This claim was wrong, 
as Theorem~\ref{th1} shows.

It is known that (the additive group of) 
the space $L^1(\mu)$ is unitarily 
representable for every measure space $(\Omega,\mu)$. (For the reader's convenience,
we remind the proof in Section~\ref{s:aux}, see Fact~\ref{f4}.) In particular, if $\mu$
is the counting measure on a set $A$, we see that the Banach space $l^1(A)$ of summable
sequences is unitarily representable. Since the product of any family of
unitarily representable groups is unitarily representable (consider the Hilbert
sum of the spaces of corresponding representations), we see that Theorem~\ref{th1}
is a consequence of the following:

\begin{theorem}
\label{th2}
For every Tikhonov space $X$ the free locally convex space $L(X)$ is isomorphic
to a subspace of a power of the Banach space $l^1(A)$ for some $A$.
\end{theorem}

For $A$ we can take any infinite set such that the cardinality of every discrete
family of non-empty open sets in $X$ does not exceed $\Card(A)$. 

Theorem~\ref{th1} implies the following result from \cite{PGao02}:
every Polish abelian group is the quotient of a closed abelian subgroup of the unitary
group of a separable Hilbert space
(the non-abelian version of the reduction is explained in Section~\ref{s:open},
the abelian case is the same). 
It is an open question (A. Kechris) whether a similar assertion holds for
non-abelian Polish groups, see \cite[Section 5.2, Question~16]{P05}, 
\cite[Question~34]{P-40Q}.

The {\em fine uniformity} $\mu_X$ on a Tikhonov space $X$ is the finest uniformity
compatible with the topology. It is generated by the family of all continuous
pseudometrics on~$X$. It is also the uniformity induced on $X$ by the group uniformity
of $A(X)$ or $L(X)$. A {\em fine uniform space} is a space of the form $(X, \mu_X)$.
We note the following corollary of Theorem~\ref{th2} which may be of some
independent interest.

\begin{corollary}\label{corol}
For every Tikhonov space $X$ the fine uniform space $(X,\mu_X)$ is isomorphic
to a uniform subspace of a power of a Hilbert space.
\end{corollary}

This need not be true for uniform spaces which are not fine. To deduce 
Corollary~\ref{corol} from Theorem~\ref{th1}, note that the left uniformity on
the unitary group $U_s(H)$ is induced by the product uniformity of $H^H$, hence
the same is true for the left uniformity on every unitarily representable group.

We prove Theorem~\ref{th2} in Section~\ref{s:proof}. The proof depends on two facts:
(1) every Banach space is a quotient of a Banach space of the form $l^1(A)$;
(2) every onto continuous map between Banach spaces admits a (possibly non-linear)
continuous right inverse. We remind the proof of these facts in Section~\ref{s:aux}. 

\section{Proof of Theorem~\ref{th2}\label{s:proof}}

Let $X$ be a Tikhonov space. Let $\sT_0$ be the topology of the free locally convex
space $L(X)$. Let $\sT_1$ be the topology on $L(X)$ generated by the linear extensions 
of all possible continuous maps of $X$ to spaces 
of the form $l^1(A)$. Theorem~\ref{th2} means that $\sT_1=\sT_0$.
In order to prove this, it suffices to verify
that $(L(X), \sT_1)$ has the following universal property: for every continuous map
$f:X\to F$, where $F$ is a Hausdorff LCS, the linear extension of $f$, say 
$\bar f:L(X)\to F$, is $\sT_1$-continuous. Since every Hausdorff LCS embeds in a product 
of Banach spaces, we may assume that $F$ is a Banach space. 
Represent $F$ as a quotient of $l^1(A)$
(Fact~\ref{f1}), let $p:l^1(A)\to F$ be linear and onto. 

$$
\xymatrix{
&l^1(A) \ar[d]_p\\
X\ar[r]_f \ar[ur]^g&F\ar@/_/@{-->}[u]_s
}
\qquad \qquad\qquad
\xymatrix{
&l^1(A) \ar[d]_p\\
L(X)\ar[r]_{\bar f} \ar[ur]^{\bar g}&F
}
$$

There exists a lift $g:X\to l^1(A)$ of $f$, i.e.{} such a map $g$ for which $f=pg$.
Indeed, find a non-linear right inverse $s:F\to l^1(A)$
of $p$ (Fact~\ref{f2}), and set $g=sf$. Then $f=psf=pg$. 
By the definition of $\sT_1$, 
the linear map $\bar g:L(X)\to l^1(A)$ is $\sT_1$-continuous. Hence $\bar f= p\bar g$
is $\sT_1$-continuous as well.

\section{Basic facts: reminders\label{s:aux}}

\begin{fact}
\label{f1}
Every Banach space $X$ is a Banach quotient of a Banach space of the form $l^1(A)$.
\end{fact}

By a Banach quotient we mean a Banach space of the form $E/F$ with the norm 
$\Vert x+F\Vert=\inf\{\Vert y\Vert: y\in x+F\}$.

\begin{proof}
Take for $A$ a dense subset of the unit ball of $X$, and consider the natural map
$p:l^1(A)\to X$. 
\end{proof}

\begin{fact}[the Bartle--Graves theorem {\cite[C.1.2]{RepSem}}]
\label{f2}

Every linear onto map $p:E\to F$ between Banach spaces (or locally convex Fr\'echet spaces)
has a (possibly non-linear) continuous 
right inverse $s:F\to E$, i.e.{} such a map that $ps=1_F$.
\end{fact}

\begin{proof}
This is a consequence of Michael's Selection Theorem for convex-valued maps 
\cite{Mi56, RepSem}. The theorem reads as follows. Suppose $X$ is paracompact, 
$E$ is a locally
convex Fr\'echet space, and for every $x\in X$ a closed convex
non-empty set $\Phi(x)\sbs E$ is given. 
Suppose further that $\Phi$ is {\em lower semicontinuous}: for every
$U$ open in $E$, the set $\Phi\obr(U)=\{x\in X: \Phi(x)\hbox{ meets }U\}$ is open in $X$.
Then $\Phi$ has a continuous selection: there exists a continuous map $s:X\to E$ such that
$s(x)\in \Phi(x)$ for every $x\in X$.

To prove Fact~\ref{f2}, apply Michael's selection theorem to $X=F$ and $\Phi$ defined
by $\Phi(x)=p\obr(x)$. The lower semicontinuity of $\Phi$ follows from the Open Mapping
theorem, according to which $p$ is open. 

For various proofs of Michael's selection theorem, see \cite{RepSem}. In particular,
note a nice reduction of the convex-valued selection theorem to the zero-dimensional
selection theorem  \cite[A, \S3]{RepSem}, based on the notion of a 
Milyutin map. 
\end{proof}

For a complex matrix $A=(a_{ij})$ we denote by $A^*$ the matrix $(\bar a_{ji})$.
$A$ is {\em Hermitian} if $A=A^*$. A Hermitian matrix $A$ is {\em positive}
if all the eigenvalues of $A$ are $\ge 0$ or, equivalently,
if $A=B^2$ for some Hermitian $B$. A complex function $p$ on a group $G$ is
{\em positive-definite}, or {\em of positive type},
if $p(g\obr)=\overline{p(g)}$ for every $g\in G$, 
and for every $g_1,\dots, g_n\in G$ the Hermitian $n\ti n$-matrix
$(p(g_i\obr g_j))$ is positive.
 If $f:G\to U(H)$ is a unitary
representation of $G$, then for every vector $v\in H$ the function $p=p_v$
on $G$ defined by $p(g)=(gv,v)$ is positive-definite. 
(We denote by $(x,y)$ the scalar product of $x,y\in H$.) 
Conversely, let $p$ be a positive-definite function on $G$.
Then $p=p_v$ for some unitary representation $f:G\to U(H)$ and some $v\in H$.
Indeed, consider the group algebra $\C[G]$, equip it with the scalar
product defined by $(g,h)=p(h\obr g)$, quotient out the kernel, and take
the completion for $H$. The regular representation of $G$ on $\C[G]$
gives rise to a unitary representation on $H$ with the required property.
(This is the so-called {\em GNS-construction}, 
see e.g.{} \cite[Theorem C.4.10]{bekka} for more details.)
If $G$ is a topological group and the function $p$ is continuous, the resulting
unitary representation is continuous as well. These arguments yield the
following (see e.g.{} \cite[Proposition 2.1]{Usp04}):

\begin{fact}
\label{f3}
A topological group $G$ is unitarily representable
(in other words, is isomorphic to a subgroup of $U_s(H)$ for some Hilbert space
$H$) if and only if for every neighbourhood $U$ of the neutral element $e$ of $G$
there exist a continuous positive-definite function $p:G\to \C$ and $a>0$
such that $p(e)=1$ and $|1-p(g)|>a$ for every $g\in G\stm U$.
\end{fact}

For a measure space $(\Omega,\mu)$ we denote by 
$L^1(\mu)$ the complex Banach space of (equivalence classes of)
complex integrable functions, and by $L_{\R}^1(\mu)$ the real Banach space
of (equivalence classes of) real integrable functions. 

\begin{fact}[Schoenberg \cite{Sch37, Sch38}]
\label{f5}
If $(\Omega,\mu)$ is a measure space and $X=L_{\R}^1(\mu)$, the function 
$x\mapsto \exp(-\Vert x\Vert)$ on $X$ is positive-definite. 
In other words, for any $f_1,\dots, f_n\in L_{\R}^1(\mu)$ the symmetric real
matrix $(\exp(-\|f_i-f_j\|))$ is positive.
\end{fact}

\begin{proof}
We invoke Bochner's theorem: positive-definite continuous functions on $\R^n$ 
(or any locally
compact abelian group) are exactly the Fourier transforms of positive measures.
For $f\in L^1(\R^n)$ we define the Fourier transform $\hat f$ by
$$
\hat f(y)=\int_{\R^n}f(x)\exp{(-2\pi i(x,y))}\,dx.
$$
Here $(x,y)=\sum_{k=1}^n x_ky_k$ for $x=(x_1,\dots,x_n)$ and  
$y=(y_1,\dots,y_n)$.

The positive functions $p$ and $q$ on $\R$ defined by $p(x)=\exp(-|x|)$ and 
$q(y)=2/(1+4\pi^2y^2)$ are the Fourier transforms of each other. Hence
each of them is positive-definite. Similarly, the positive functions $p_n$ and $q_n$
on $\R^n$ defined by $p_n(x_1,\dots,x_n)=\exp(-\sum_{k=1}^n |x_k|)=
\prod_{k=1}^n p(x_k)$ and 
$q_n(y_1,\dots,y_n)=\prod_{k=1}^n q(y_k)$ are positive-definite, being
the Fourier transforms of each other. If $m_1,\dots,m_n$ are strictly positive
masses, 
the function $x\mapsto \exp(-\sum_{k=1}^n m_k|x_k|)$ is positive-definite, since it is
the composition of $p_n$ and a linear automorphism of $\R^n$. This is exactly
Fact~\ref{f5} for finite measure spaces.

The general case easily follows: given finitely many functions 
$f_1,\dots, f_n\in L_{\R}^1(\mu)$, we can approximate them by finite-valued functions.
In this way we see 
that the symmetric matrix $A=(\exp(-\|f_i-f_j\|))$ is in the closure of the set of matrices
$A'$ of the same form arising from finite measure spaces. The result of the preceding
paragraph means that each 
$A'$ is positive. Hence $A$ is positive.
\end{proof}

As a topological group, $L^1(\mu)$ is isomorphic to the square of $L_{\R}^1(\mu)$.
Combining Facts~\ref{f3} and~\ref{f5},
we obtain:

\begin{fact}
\label{f4}
The additive group of
the space $L^1(\mu)$ is unitarily 
representable for every measure space $(\Omega,\mu)$.
\end{fact}

See \cite{Megr} for more on unitarily and reflexively representable Banach spaces.

\section{Open questions\label{s:open}}

Let us say that a metric space $M$ is of {\em $L^1$-type} if it is isometric to 
a subspace of the Banach space $L_{\R}^1(\mu)$ for some measure space $(\Omega, \mu)$. 
A non-abelian version of Theorems~\ref{th1} and~\ref{th2} might be the following:

\begin{conjecture}\label{conj}
For any Tikhonov space $X$ the free topological group $F(X)$ is isomorphic to a subgroup
of the group of isometries $\Iso(M)$ for some metric space $M$
of $L^1$-type.
\end{conjecture}

It follows from \cite[Theorem 3.1]{Usp04} (apply it to the positive-definite
function $p$ on $\R$ used in the proof of Fact~\ref{f5} and
defined by $p(x)=\exp(-|x|)$) and Fact~\ref{f5} 
that for every $M\sbs L_{\R}^1(\mu)$ the group $\Iso(M)$
is unitarily representable. Thus, if conjecture~\ref{conj} is true,
every $F(X)$ is unitarily representable. This would imply a positive answer to 
the question of Kechris mentioned in Section~\ref{s:intro}: is
every Polish group a quotient
of a closed subgroup of the unitary group of a separable Hilbert space?
Indeed:

\begin{proposition}\label{prop}
Let $P$ be the space of irrationals. If 
the group $F(P)$ is unitarily representable, then every Polish group is a quotient
of a closed subgroup of the unitary group of a separable Hilbert space.
\end{proposition}

\begin{proof}
A topological group is {\em uniformly Lindel\"of} 
(or, in another terminology, {\em $\o$-bounded}) if for every neighbourhood $U$ 
of the unity the group can be covered by countably many left (equivalently, right)
translates of $U$. If $G$ is a uniformly Lindel\"of group of isometries of a metric
space $M$, then for every $x\in M$ the orbit $Gx$ is separable
(see e.g.{} the section ``Guran's theorems" in \cite{Usp86}).
If $G$ is a uniformly Lindel\"of subgroup of the unitary group $U_s(H)$, where
$H$ is a (non-separable) Hilbert space, it easily follows that $H$ is covered by 
separable closed $G$-invariant linear subspaces and therefore $G$ 
embeds in a product of unitary groups of separable Hilbert spaces.

The group $F(P)$, like any separable topological group, is uniformly Lindel\"of.
If we assume that $F(P)$ is unitarily representable, it follows that $F(P)$
is isomorphic to a topological subgroup
of a power of $U_s(H)$, where $H$ is a separable Hilbert space.

If $G$ is a Polish group,
there exists a quotient onto map $F(P)\to G$ (because there exists an open onto map
$P\to G$), and an easy factorization argument shows that there is a group $N$ lying in
a countable power of $U_s(H)$ (and hence isomorphic to a subgroup of $U_s(H)$) such that
$G$ is a quotient of $N$. The quotient homomorphism $N\to G$ can be extended over 
the closure of $N$, so we may assume that $N$ is closed in $U_s(H)$.
\end{proof}

For any uniform space $X$ one defines the free abelian group $A(X)$ and free locally
convex space $L(X)$ in an obvious way. The objects $A(X)$ and $L(X)$ for 
Tikhonov spaces $X$ considered in this paper are special cases of the same objects
for uniform spaces, 
corresponding to fine uniform spaces.

\begin{question}[Megrelishvili]
\label{que}
For what uniform spaces $X$ are the groups $A(X)$ and $L(X)$ unitarily representable?
Is it sufficient that $X$ be a uniform subspace of a product of Hilbert spaces?
\end{question}

\section{Acknowledgement}
I thank my friends M. Megrelishvili and V. Pestov for their generous help,
which included inspiring discussions and many
useful remarks and suggestions on early versions of this paper.

\end{document}